\documentclass[11pt]{article}
\usepackage{amssymb,amsmath,bm,graphicx}
\usepackage{mathrsfs}
\usepackage{constants}
\usepackage{hyperref}
\newcommand{\bsxi}{\boldsymbol{\xi}}

\usepackage{color}

\date{}
\begin{document}
\newcommand{\bea}{\begin{eqnarray}}
\newcommand{\ena}{\end{eqnarray}}
\newcommand{\beas}{\begin{eqnarray*}}
\newcommand{\enas}{\end{eqnarray*}}
\newcommand{\beq}{\begin{equation}}
\newcommand{\enq}{\end{equation}}
\def\qed{\hfill \mbox{\rule{0.5em}{0.5em}}}
\newcommand{\bbox}{\hfill $\Box$}
\newcommand{\ignore}[1]{}
\newcommand{\ignorex}[1]{#1}
\newcommand{\wtilde}[1]{\widetilde{#1}}
\newcommand{\qmq}[1]{\quad\mbox{#1}\quad}
\newcommand{\qm}[1]{\quad\mbox{#1}}
\newcommand{\nn}{\nonumber}
\newcommand{\Bvert}{\left\vert\vphantom{\frac{1}{1}}\right.}
\newcommand{\To}{\rightarrow}
\newcommand{\e}{\mathbb{E}}
\newcommand{\Var}{\mathrm{Var}}
\newcommand{\Cov}{\mathrm{Cov}}
\makeatletter
\newsavebox\myboxA
\newsavebox\myboxB
\newlength\mylenA
\newcommand*\xoverline[2][0.70]{%
    \sbox{\myboxA}{$\m@th#2$}%
    \setbox\myboxB\null
    \ht\myboxB=\ht\myboxA%
    \dp\myboxB=\dp\myboxA%
    \wd\myboxB=#1\wd\myboxA
    \sbox\myboxB{$\m@th\overline{\copy\myboxB}$}
    \setlength\mylenA{\the\wd\myboxA}
    \addtolength\mylenA{-\the\wd\myboxB}%
    \ifdim\wd\myboxB<\wd\myboxA%
       \rlap{\hskip 0.5\mylenA\usebox\myboxB}{\usebox\myboxA}%
    \else
        \hskip -0.5\mylenA\rlap{\usebox\myboxA}{\hskip 0.5\mylenA\usebox\myboxB}%
    \fi}
\makeatother

\newtheorem{theorem}{Theorem}[section]
\newtheorem{corollary}[theorem]{Corollary}
\newtheorem{conjecture}[theorem]{Conjecture}
\newtheorem{proposition}[theorem]{Proposition}
\newtheorem{lemma}[theorem]{Lemma}
\newtheorem{definition}[theorem]{Definition}
\newtheorem{example}[theorem]{Example}
\newtheorem{remark}[theorem]{Remark}
\newtheorem{case}{Case}[section]
\newtheorem{condition}{Condition}[section]
\newcommand{\proof}{\noindent {\it Proof:} }

\title{{\bf\Large Stein's method for negatively associated random variables with applications to second order stationary random fields}}
\author{Nathakhun Wiroonsri \\University of Southern California}
\footnotetext{Department of Mathematics, University of Southern California, Los Angeles, CA 90089, USA, wiroonsr@usc.edu}
\footnotetext{Research supported by the 2017-2018 Russell Endowed Fellowship, University of Southern California.}
\maketitle

\begin{abstract}   
Let $\bsxi=(\xi_1,\ldots,\xi_m)$ be a negatively associated mean zero random vector with components that obey the bound $|\xi_i| \le B, i=1,\ldots,m$, and whose sum $W = \sum_{i=1}^m \xi_i$ has variance 1, the bound
\beas
d_1\big({\cal L}(W),{\cal L}(Z)\big) \le  5B - 5.2\sum_{i \not = j} \sigma_{ij}.
\enas
is obtained where $Z$ has the standard normal distribution and $d_1(\cdot,\cdot)$ is the $L^1$ metric. The result is extended to the multidimensional case with the $L^1$ metric replaced by a smooth functions metric. Applications to second order stationary random fields with exponential decreasing covariance are also presented.
\end{abstract}

\section{Introduction}

This work is extended from the recent paper \cite{GW16} that focuses only on positive association. Particularly, in this work, we provide non-asymptotic $L^1$ bounds to the normal for negatively associated random fields in $\mathbb{Z}^d$ using the same technique developed in \cite{GW16}. We recall that the $L^1$, or Wasserstein, distance between the distributions ${\cal L}(X)$ and ${\cal L}(Y)$ of real valued random variables $X$ and $Y$ is given by
\bea \label{def.d1.integral}
d_1 \big({\cal L}(X),{\cal L}(Y)\big) = \int_{-\infty}^\infty |P(X \le t) - P(Y \le t)| dt.
\ena 

A random vector $\bsxi=(\xi_1,\ldots,\xi_m) \in \mathbb{R}^m$ is said to be \textit{positively associated} whenever
\beas
\Cov\big(f(\bsxi),g(\bsxi)\big) \ge 0
\enas
for all real valued coordinate-wise nondecreasing functions $f$ and $g$ on $\mathbb{R}^m$ such that 
$f(\bsxi)$ and $g(\bsxi)$ have finite second moments. A random vector $\bsxi=(\xi_1,\ldots,\xi_m) \in \mathbb{R}^m$ is said to be \textit{negatively associated} if for any disjoint subsets $A,B$ of $[m] := \{1,2,\ldots,m \}$,
\bea \label{negdef}
\Cov\big(f(\xi_i,i\in A),g(\xi_j,j\in B)\big) \le 0
\ena
for all coordinate-wise increasing functions $f:\mathbb{R}^{|A|} \rightarrow \mathbb{R}$ and $g:\mathbb{R}^{|B|} \rightarrow \mathbb{R}$ such that $f(\xi_i,i\in A)$ and $g(\xi_j,j\in B)$ have finite second moments. In general, a collection $\{\xi_\alpha: \alpha \in I\}$ of real valued random variables indexed by a set $I$ is said to be positively associated (resp. negatively associated) if all finite subcollections are positively associated (resp. negatively associated). Positive association was introduced in \cite{EPW67} and has been found frequently in probabilistic models in several areas, especially statistical physics. In some literature positive association is termed the `FKG-inequality' or simply `association' (see \cite{new80} and \cite{cox84} for examples). Negative association was later introduced in \cite{JP83} and has well-known applications related to permutation distributions. 

Over the last few decades, many researchers established central limit theorems and rates of convergence for sums of positively associated random variables (\cite{new80},\cite{Bir88},\cite{Bul95}) and negatively associated random variables (\cite{LW08},\cite{CW09}) under different assumptions. Recently, the work \cite{GW16} developed an $L^1$ version of Stein's method adapted to sums of positively associated random variables with applications to statistical physics. Stein's method, introduced by \cite{Stein72}, is nowadays one of the most powerful methods to prove convergence in distribution as it has main advantages that it provides non-asymptotic bounds on the distance between distributions, and that it can handle various situations involving dependence. Thus far, many applications in several areas such as statistics, statistical physics and applied sciences have been developed using this method. For more detail about the method in general, see the text \cite{Che11} and the introductory notes \cite{Ross11}. 

The one dimensional result of \cite{GW16} is stated below. 

\begin{theorem}[\cite{GW16}] \label{gw1}
Let $\bsxi=(\xi_1,\ldots,\xi_m)$ be a positively associated mean zero random vector with components obeying the bound $|\xi_i| \le B$ for some $B>0$, and whose sum $W = \sum_{i=1}^m \xi_i$ has variance 1. Let $Z$ be a standard normal random variable. Then
\beas 
d_1\big({\cal L}(W),{\cal L}(Z)\big) \le  5B + \sqrt{\frac{8}{\pi}}\sum_{i \not = j} \sigma_{ij} \qmq{where} \sigma_{ij}=\e[\xi_i \xi_j].
\enas
\end{theorem}

The multidimensional result was also obtained in \cite{GW16} with the $L^1$ metric replaced by a smooth functions metric, following the development of Chapter 12 of \cite{Che11}. 

For 
\beas
{\bf x} \in \mathbb{R}^p \qmq{let} |{\bf x}|_1=\sum_{i=1}^p |x_i|, \quad \mbox{the $L^1$ vector norm,}
\enas
and for a real valued function $\varphi(u)$ defined on the domain $\mathcal{D}$, let $|\varphi|_{\infty}= \sup_{x \in \mathcal{D}} |\varphi(x)|$. We include in this definition the $|\cdot|_\infty$ norm of vectors and matrices, for instance, by considering them as real valued functions of their indices. Also, from this point, we denote $\mathbb{N}_k=[k,\infty) \cap \mathbb{Z}$ for $k \in \mathbb{Z}$.

For $m \in \mathbb{N}_0$, let $L^{\infty}_m(\mathbb{R}^p)$ be the collection of all functions $h:\mathbb{R}^p \rightarrow \mathbb{R}$ such that for all $\mathbf{k} = (k_1,\ldots,k_p) \in \mathbb{N}_0^p$ with $|{\bf k}|_1 \le m$, the partial derivative 
\beas
h^{(\mathbf{k})}(\mathbf{x}) = \frac{\partial^{|\mathbf{k}|_1}h}{\partial^{k_1}x_1\cdots\partial^{k_p}x_p} 
\enas
exists, and 
\beas
|h|_{L^{\infty}_m(\mathbb{R}^p)}: = \max_{0 \le |\mathbf{k}|_1 \le m}|h^{(\mathbf{k})}|_{\infty} \text{ \ \ is finite.}
\enas
For $f \in L_m^\infty(\mathbb{R}^p)$ let 
\beas
\mathcal{H}_{m,\infty,p} = \left\{ h\in L^{\infty}_m(\mathbb{R}^p):|h|_{L^{\infty}_m(\mathbb{R}^p)} \le 1  \right\},
\enas 
and for random vectors $\mathbf{X}$ and $\mathbf{Y}$ in $\mathbb{R}^p$, 
define the smooth functions metric
\bea \label{smoothdef}
d_{\mathcal{H}_{m,\infty,p}}\big(\mathcal{L}(\mathbf{X}),\mathcal{L}(\mathbf{Y})\big) = \sup_{h \in {{\cal H}}_{m,\infty,p}} |\e h(\mathbf{X})-\e h(\mathbf{Y})|.
\ena

For a positive semidefinite matrix $H$, we let $H^{1/2}$ denote the unique positive semidefinite square root of $H$. When $H$ is positive definite, we write $H^{-1/2}=(H^{1/2})^{-1}$. The following theorem states the multidimensional result of \cite{GW16}. 

\begin{theorem}[\cite{GW16}]  \label{gw2}
With $m,p \in \mathbb{N}_1$, let $\{ \xi_{i,j}:i \in [m],j\in [p] \}$ be positively associated mean zero random variables bounded in absolute value by some positive constant $B$. Let $\mathbf{S} = (S_1,S_2,\ldots,S_p)$ where $S_j= \sum_{1 \le i \le m}\xi_{i,j}$ for $j \in [p]$ and assume that $\Sigma = \Var \big({\bf S} \big)$ is positive definite. Then
\beas
&&d_{\mathcal{H}_{3,\infty,p}}\big(\mathcal{L}(\Sigma^{-1/2}\left( \mathbf{S}-\e \mathbf{S})\right),\mathcal{L}(\mathbf{Z})\big) \\
&& \hspace{20pt}\le  \left(\frac{1}{6}+2\sqrt{2}\right)p^3B|\Sigma^{-1/2}|_{\infty}^3\sum_{j=1}^p\Sigma_{j,j} \\
&&\hspace{30pt}+ \left(\frac{3}{\sqrt{2}}+\frac{1}{2}\right)p^2|\Sigma^{-1/2}|_{\infty}^2 \sum_{j=1}^p \sum_{i,k\in[m],i \ne k} \Cov\left(\xi_{i,j},\xi_{k,j}\right) \\
&&\hspace{30pt}+ \left( 2\sqrt{2}p^3B |\Sigma^{-1/2}|_{\infty}^3 + \left(\frac{3}{\sqrt{2}}+\frac{1}{2}\right)p^2|\Sigma^{-1/2}|_{\infty}^2 \right) \sum_{j,l \in [p],j \ne l} \Sigma_{j,l},
\enas
where $\mathbf{Z} \sim {\cal N}({\bf 0},I_p)$, a standard normal vector in $\mathbb{R}^p$.
\end{theorem}

The results above were applied to second order stationary random fields assuming exponential decreasing covariance and to four models in statistical physics; Ising and voter models, bond percolation and contact process. In the present work, we prove similar results adapted to negative association. 

Stein's method has been used previously for negative association and some related concepts. In \cite{LW08}, Stein's method was used in normal approximation for sums of pairwise negative quadrant dependent random variables which allows one to derive a CLT for pairwise negative quadrant dependent random variables with Lindeberg's condition. We note that when $m=2$ in \eqref{negdef} negative association and negative quadrant dependent are equivalent (see \cite{JP83}). In \cite{Dal13}, Stein's method was used to obtain the total variation distance between compound Poisson distribution and sums of positively associated or negatively associated random variables. 

Our main results in the one dimensional and multidimensional cases are stated in Theorems \ref{stein:fkg} and \ref{stein:fkg2}, respectively, as follows.

\begin{theorem}
\label{stein:fkg}
Let $\bsxi=(\xi_1,\ldots,\xi_m)$ be a negatively associated mean zero random vector with components obeying the bound $|\xi_i| \le B$ for some $B>0$, and whose sum $W = \sum_{i=1}^m \xi_i$ has variance 1. Let $Z$ be a standard normal random variable. Then, with $\sigma_{ij}=\e[\xi_i \xi_j]$,
\bea \label{stein:fkgbound}
d_1\big({\cal L}(W),{\cal L}(Z)\big) \le  5B - 5.2\sum_{i \not = j} \sigma_{ij}.
\ena
\end{theorem}

\begin{theorem} \label{stein:fkg2}
With $m,p \in \mathbb{N}_1$, let $\{ \xi_{i,j}:i \in [m],j\in [p] \}$ be negatively associated mean zero random variables satisfying $|\xi_{i,j}|\le B$ for some $B>0$. Let $\mathbf{S} = (S_1,S_2,\ldots,S_p)$ where $S_j= \sum_{1 \le i \le m}\xi_{i,j}$ for $j \in [p]$ and assume that $\Sigma = \Var \big({\bf S} \big)$ is positive definite. Then
\bea\label{steinbound:fkg2}
&&d_{\mathcal{H}_{3,\infty,p}}\big(\mathcal{L}(\Sigma^{-1/2}\left( \mathbf{S}-\e \mathbf{S})\right),\mathcal{L}(\mathbf{Z})\big) \nn\\
&&\hspace{30pt}\le  \frac{5}{6}p^3B|\Sigma^{-1/2}|_{\infty}^3\sum_{j=1}^p\Sigma_{j,j}  \nn \\
&&\hspace{40pt}- \left(\frac{3}{2}p^3B|\Sigma^{-1/2}|_{\infty}^3 + p^2|\Sigma^{-1/2}|_{\infty}^2\right) \sum_{j=1}^p \sum_{i,k\in[m],i \ne k} \Cov\left(\xi_{i,j},\xi_{k,j}\right) \nn \\
&&\hspace{40pt}- \left( \frac{2}{3}p^3B |\Sigma^{-1/2}|_{\infty}^3 + p^2|\Sigma^{-1/2}|_{\infty}^2 \right) \sum_{j,l \in [p],j \ne l} \Sigma_{j,l} ,
\ena
where $\mathbf{Z} \sim {\cal N}({\bf 0},I_p)$, a standard normal random vector in $\mathbb{R}^p$.
\end{theorem}

\begin{remark}
We note that the differences between our results in Theorems \ref{stein:fkg} and \ref{stein:fkg2} for negative association and those in Theorems \ref{gw1} and \ref{gw2} of \cite{GW16} for positive association, are that the signs of the covariance terms are reverse and that the constants are different. Nevertheless, in Section \ref{app} we have an example where these changes do not contribute rates of convergence. We also note that the bounds in the four theorems above are particularly useful when the variables one handles are bounded and (positively or negatively) associated. However let us compare the one dimensional results in Theorems \ref{gw1} and \ref{stein:fkg} with the classical result for independent and identically distributed variables $X_i$, $i\in [n]$ with $\e X_i =0$, $\Var X_i = \sigma^2>0$ and $|X_i| \le K$ for some $K>0$. With $W = \left(\sum_i X_i\right)/(\sigma\sqrt{n})$, both Theorems \ref{gw1} and \ref{stein:fkg} give that $d_1\big({\cal L}(W),{\cal L}(Z)\big) \le 5K\sigma^{-1}n^{-1/2}$. The classical Berry Esseen theorem provides the bound $\sup_{x \in \mathbb{R}} |P(W \le x) - P(Z \le x)| \le CK^3\sigma^{-3}n^{-1/2}$ with the smallest constant $C = 0.4748$ obtained recently by \cite{She11}. Though the rates of convergence in $n$ are the same, the constants are different and actually the two distances are not comparable as it is known that $\sup_{x \in \mathbb{R}} |P(X \le x) - P(Y \le x)| \le \sqrt{2cd_1\big({\cal L}(X),{\cal L}(Y)\big)}$ for some $c>0$ but not conversely (See Proposition 1.2 of \cite{Ross11}). 
\end{remark}

The remainder of this work is organized as follows. In the next section, we use Stein's method to prove the two main results, Theorems \ref{stein:fkg} and \ref{stein:fkg2}. We state our results for negatively associated random fields whose covariance decays exponentially in Section \ref{app} . One similar advantage of the four theorems stated in this section is that, unlike many results based on Stein's method, they may be applied without the need for coupling constructions.

\section{Proofs of main theorems} 

In this section we prove our main results, Theorems \ref{stein:fkg} and \ref{stein:fkg2}, using similar techniques as in \cite{GW16}. For this purpose, we first state the following two lemmas proved in \cite{JP83} and \cite{CW09}, respectively. The more general version of the first lemma was originally proved in \cite{BS98} in Russian and the English version can be found in the book \cite{BS07}. Lemma \ref{lemma:inc} states that two disjoint sums of negatively associated random variables are negatively associated which will be used throughout the remainder of this work. Lemma \ref{lemma:covbound} allows us to bound $\left|\Cov(f(\xi_i:i\in A),g(\xi_j:j\in B))\right|$ by a linear combination of $-\Cov(\xi_i,\xi_j)$ with $i \in A, j \in B$ when $A$ and $B$ are disjoint.

\begin{lemma}[\cite{JP83}] \label{lemma:inc} 
Increasing functions defined on disjoint subsets of a set of negatively associated random variables are negatively associated.
\end{lemma}

\begin{lemma}[\cite{BS98}, \cite{CW09}] \label{lemma:covbound} 
Let $A$ and $B$ be disjoint finite sets, and let $\xi_j$, $j \in A \cup B$, be negatively associated random variables. If $f:\mathbb{R}^{|A|} \rightarrow \mathbb{R}$, $g:\mathbb{R}^{|B|} \rightarrow \mathbb{R}$ are partially differentiable with bounded  partial derivatives, then 
\beas
\left|\Cov(f(\xi_i:i\in A),g(\xi_j:j\in B))\right| \le -\sum_{i \in A} \sum_{j \in B}\left|\frac{\partial f}{\partial \xi_i}\right|_{\infty} \left|\frac{\partial g}{\partial \xi_j}\right|_{\infty} \Cov(\xi_i,\xi_j).
\enas
\end{lemma}

We note that the difference between the proofs below and the ones in \cite{GW16} results from that Lemma \ref{lemma:covbound} requires $A$ and $B$ to be disjoint unlike the one for positive association. Therefore we add a few more steps in the proofs to handle this situation. In the proof that follows, we use the alternate form of the $L^1$, or Wasserstein distance (see e.g. \cite{Rac84});
\bea \label{d1:sup.L}
d_1\big({\cal L}(X),{\cal L}(Y)\big) = \sup_{h \in \mathbb{L}} |\e h(X)-\e h(Y)| \qmq{where} \mathbb{L}=\{h: |h(y)-h(x)| \le |y-x|\}.
\ena

\noindent {\bf Proof of Theorem \ref{stein:fkg}}  For given $h \in \mathbb{L}$ let $f$ be the unique bounded solution to the Stein equation
\bea \label{stein:eq}
f'(w)-w f(w) = h(w)-Nh \qmq{where} N h=\e h(Z),
\ena
with ${\cal L}(Z)$ the standard normal distribution. Then, (see e.g. Lemma 2.4 of \cite{Che11}),
\bea \label{steineq:bounds}
|f'|_{\infty} \le \sqrt{\frac{2}{\pi}} \qmq{and} |f''|_{\infty} \le 2.
\ena

Recall that in the proof below we use the notations $\sigma_{ij}=\e[\xi_i \xi_j]$ and $\sigma^2_{i}=\Var(\xi_i)$ for $i\ne j \in [n]$. As ${\rm Var}(W)=\sum_{i=1}^n\sigma_i^2 + \sum_{i\ne j} \sigma_{ij}=1$, we obtain
\beas
\e[f'(W)]&=& \e\left(\sum_{i=1}^m \sigma_i^2 f'(W) + \sum_{i \not = j} \sigma_{ij} f'(W) \right)\\
         &=& \e\left(\sum_{i=1}^m \xi_i^2 f'(W)+ \sum_{i \not = j} \sigma_{ij} f'(W) + \sum_{i=1}^m (\sigma_i^2-\xi_i^2) f'(W) \right).
\enas 
Now letting $W^i=W-\xi_i$, write
\beas
\e[Wf(W)] = \e\sum_{i=1}^m \xi_i f(W) = \e\sum_{i=1}^m \xi_i f(W^i  + \xi_i)
         = \e\sum_{i=1}^m \left[ \xi_i f(W^i) + \xi_i^2 \int_0^1 f'(W^i + u \xi_i) du \right].
\enas
Recalling the Stein equation (\ref{stein:eq}) and subtracting two equations above, we obtain
\begin{bea} \label{add.subtract.alt}
&&\e[h(W)-N h]=\e[f'(W)-Wf(W)]\nn\\
&&\hspace{50pt}= \e \left( \sum_{i=1}^m \xi_i^2 \left(  \int_0^1 \left(f'(W) -f'(W^i+ u\xi_i)\right)du\right)+  \sum_{i=1}^m (\sigma_i^2 - \xi_i^2) f'(W) \right.\nn\\
&&\hspace{100pt}\left.  + \sum_{i \not = j} \sigma_{ij}f'(W) - \sum_{i=1}^m  \xi_i f(W^i)\right).
\ena

Using the second inequality in (\ref{steineq:bounds}),  we bound the first term in (\ref{add.subtract.alt}) by
\bea \label{term:1}
&&\left| \e \sum_{i=1}^m \xi_i^2 \int_0^1 \left(f'(W) -f'(W^i+ u\xi_i) \right) du \right| \nn\\
&&\hspace{50pt}= \left| \e \sum_{i=1}^m \xi_i^2 \int_0^1 \int_{u\xi_i}^{\xi_i} f''(W^i+t) dt du  \right|
\le 2 \e \sum_{i=1}^m \xi_i^2 \left(  \int_0^1 \int_{u|\xi_i|}^{|\xi_i|} dt du\right) \nn\\
&&\hspace{50pt}=  \e \sum_{i=1}^m |\xi_i|^3 \le B \e \sum_{i=1}^m \xi_i^2 = B\left(1-\sum_{i \ne j} \sigma_{ij}\right).
\ena

To handle the second term in (\ref{add.subtract.alt}), using the triangle inequality, we first bound it by the three terms denoted by $I_1,I_2,I_3$, respectively, 
\beas 
\left| \e \sum_{i=1}^m f'(W)(\sigma_i^2 - \xi_i^2) \right|
&\le& \left|\e \sum_{i=1}^m f'(W^i)(\sigma_i^2 - \xi_i^2) \right|+\left|\e \sum_{i=1}^m (f'(W)-f'(W^i))(\sigma_i^2 - \xi_i^2) \right|\\
&\le& \left|\e \sum_{i=1}^m f'(W^i)(\sigma_i^2 - \xi_i^2) \right|+\e \sum_{i=1}^m \left|f'(W)-f'(W^i)\right|\left|\sigma_i^2-\xi_i^2\right| \\
&\le& \left|\e \sum_{i=1}^m f'(W^i)(\sigma_i^2 - \xi_i^2) \right|+\e \sum_{i=1}^m \left|f'(W)-f'(W^i)\right|\sigma_i^2 \\
    &&+\e \sum_{i=1}^m \left|f'(W)-f'(W^i)\right|\xi_i^2  := I_1+I_2+I_3.
\enas
Note that $W^i$ and $\xi_i$ are coordinate-wise increasing functions defined on disjoint subsets of $\bsxi$, and hence negatively associated by Lemma \ref{lemma:inc}. Now for $I_1$, applying Lemma \ref{lemma:covbound}  with
\beas
g(x)=\left\{  
\begin{array}{cc}
x^2 & |x| \le B\\
B^2 & |x|>B
\end{array}
\right.
\enas
and using the second inequality in (\ref{steineq:bounds}), we have
\bea
I_1 = \left| \sum_{i=1}^m \Cov \big(f'(W^i),g(\xi_i) \big) \right| \le -4B \sum_{i=1}^m   \Cov \big(W^i,\xi_i \big) = -4B \sum_{i\ne j} \sigma_{ij}.
\ena
For $I_2$ and $I_3$, applying again the second inequality in (\ref{steineq:bounds}), we obtain
\bea
I_2
\le |f''|_{\infty}\sum_{i=1}^m \sigma_i^2 \e|W-W^i| 
\le 2B\sum_{i=1}^m \sigma_i^2 =2B\left(1-\sum_{i \ne j} \sigma_{ij}\right),
\ena
and
\bea
I_3 \le  |f''|_{\infty} \sum_{i=1}^m  \e\xi_i^2|W-W^i|  \le 2B \sum_{i=1}^m  \e\xi_i^2 = 2B\left(1-\sum_{i \ne j} \sigma_{ij}\right).
\ena

For the third term in (\ref{add.subtract.alt}), using the negativity of the covariances $\sigma_{ij}$, $i\ne j$, and the first inequality in (\ref{steineq:bounds}) we obtain
\bea 
\left| \e\sum_{i \not = j} \sigma_{ij}f'(W) \right| \le -|\e f'(W)| \sum_{i \not = j} \sigma_{ij}  \le -\sqrt{\frac{2}{\pi}} \sum_{i \not = j} \sigma_{ij}.
\ena

For the final term in (\ref{add.subtract.alt}), using again the fact that the pair $(W^i,\xi_i)$ is negatively associated and applying Lemma \ref{lemma:covbound} and the first inequality in (\ref{steineq:bounds}) now yields
\bea \label{term:4}
\left| \e \sum_{i=1}^m  \xi_i f(W^i) \right| = \left| \sum_{i=1}^m \Cov\big(\xi_i,f(W^i) \big) \right| \le -\sqrt{\frac{2}{\pi}}\sum_{i=1}^m \Cov \big(\xi_i,W^i \big) = -\sqrt{\frac{2}{\pi}}\sum_{i \not = j}\sigma_{ij}.
\ena

Summing the bounds (\ref{term:1})-(\ref{term:4}), taking supremum over $h \in \mathbb{L}$ and using the form of the $L^1$ distance given in (\ref{d1:sup.L}), we obtain
\beas
d_1\big({\cal L}(W),{\cal L}(Z)\big) \le  5B - 9B\sum_{i \not = j} \sigma_{ij}-\sqrt{\frac{8}{\pi}}\sum_{i \not = j} \sigma_{ij}.
\enas
Using the fact that $d_1(\cdot,\cdot) \le 2$ and that $\sigma_{i,j}$ are negative, we can assume that $B \le 0.4$ and thus the last expression is bounded by the right hand side of (\ref{stein:fkgbound}). 

\bbox

Next we use the following result which is slightly different from Lemma 2.6 of \cite{Che11} due to \cite{Bar90} to prove Theorem \ref{stein:fkg2}. Let $\mathbf{Z}$ be a standard normal random vector in $\mathbb{R}^p$. For $h: \mathbb{R}^p \rightarrow \mathbb{R}$ let $N h = \e h(\mathbf{Z})$ and for $u \ge 0$ define
\beas
(T_uh)(\mathbf{s})=\e{h(\mathbf{s}e^{-u}+\sqrt{1-e^{-2u}}\mathbf{Z})}.
\enas
We write $D^2 h$ for the Hessian matrix of $h$ when it exists. 

\begin{lemma} \label{multilemma}
For $m \ge 3$ and $h \in L^{\infty}_m(\mathbb{R}^p)$ the function 
\beas
g(\mathbf{s}) = -\int_0^{\infty}[T_uh(\mathbf{s})-N h]du
\enas
solves
\beas
{\rm tr} D^2g(\mathbf{s})-\mathbf{s}\cdot \nabla g(\mathbf{s}) = h(\mathbf{s})-N h,
\enas
and for any $0 \le |{\bf k}|_1 \le m$
\beas
|g^{(\mathbf{k})}|_{\infty}\le \frac{1}{|{\bf k}|_1} |h^{(\mathbf{k})}|_{\infty}.
\enas
Furthermore, for any $\boldsymbol\lambda \in \mathbb{R}^p$ and positive definite $p \times p$ matrix $\Sigma$, $f$ defined by the change of variable
\bea \label{solstein}
f(\mathbf{s}) = g(\Sigma^{-1/2}(\mathbf{s}-\boldsymbol\lambda))
\ena
solves
\bea \label{Stein's eq2}
tr\Sigma D^2 f(\mathbf{s})-(\mathbf{s}-\boldsymbol\lambda) \cdot \nabla f(\mathbf{s}) = h(\Sigma^{-1/2}(\mathbf{s}- \boldsymbol\lambda)) - N h,
\ena
and satisfies
\beas
|f^{(\mathbf{k})}|_{\infty} \le \frac{p^{|\mathbf{k}|_1}}{|\mathbf{k}|_1}|\Sigma^{-1/2}|_{\infty}^{|\mathbf{k}|_1}|h^{(\mathbf{k})}|_{\infty}.
\enas
In particular, if $h \in \mathcal{H}_{m,\infty,p}$ then 
\bea \label{diffbound}
|f^{(\mathbf{k})}|_{\infty} \le \frac{p^{|\mathbf{k}|_1}}{|\mathbf{k}|_1}|\Sigma^{-1/2}|_{\infty}^{|\mathbf{k}|_1} \text{ \ for all \ } 0 \le |\mathbf{k}|_1 \le m.
\ena
\end{lemma}

We apply the same technique as in the univariate case, along with Lemmas \ref{lemma:covbound} and \ref{multilemma}, to prove our main multivariate theorem below.

\noindent {\bf Proof of Theorem \ref{stein:fkg2}} 
Given $h \in \mathcal{H}_{3,\infty,p}$, let $f$ be the solution of (\ref{Stein's eq2}) given by (\ref{solstein}) with $\boldsymbol\lambda={\bf 0}$. Writing out the expressions in (\ref{Stein's eq2}) yields
\begin{multline} \label{1st}
\e\left[h(\Sigma^{-1/2}\mathbf{S})- N h \right] = \e\left[\sum_{j=1}^p \sum_{l=1}^p \Sigma_{j,l}\frac{\partial ^2}{\partial s_j\partial s_l}f(\mathbf{S})-\sum_{j=1}^pS_j\frac{\partial}{\partial s_j}f(\mathbf{S})\right]  \\
= \e\sum_{j=1}^p\Sigma_{j,j}\frac{\partial ^2}{\partial s_j^2}f(\mathbf{S})+ \e\sum_{j,l \in [p],j \ne l} \Sigma_{j,l}\frac{\partial ^2}{\partial s_j\partial s_l}f(\mathbf{S})-\e\sum_{j=1}^pS_j\frac{\partial}{\partial s_j}f(\mathbf{S}).
\end{multline}
We consider the first term of (\ref{1st}) and deal with each term under the sum separately for $j=1, \ldots,p$. Letting $\sigma^2_{i,j}:=\Var\big(\xi_{i,j} \big)$ and $\sigma_{i,j;k,l} := \Cov\big(\xi_{i,j},\xi_{k,l}\big)$, we have
\bea \label{2nd}
\Sigma_{j,j}\frac{\partial ^2}{\partial s_j^2}f(\mathbf{S}) &=& \sum_{i=1}^m \sigma_{i,j}^2\frac{\partial ^2}{\partial s_j^2}f(\mathbf{S}) + \sum_{i,k \in [m],i \neq k} \sigma_{i,j;k,j}\frac{\partial ^2}{\partial s_j^2}f(\mathbf{S}) \nn\\ 
&=& \sum_{i=1}^m \xi_{i,j}^2\frac{\partial ^2}{\partial s_j^2}f(\mathbf{S}) + \sum_{i,k \in [m],i \neq k} \sigma_{i,j;k,j}\frac{\partial ^2}{\partial s_j^2}f(\mathbf{S}) +\sum_{i=1}^m (\sigma_{i,j}^2-\xi_{i,j}^2)\frac{\partial ^2}{\partial s_j^2}f(\mathbf{S}). 
\ena
Now, with $S_{j*i} := S_j-\xi_{i,j}$ we write the summands of the third term on the right hand side of (\ref{1st}) as
\bea \label{3rd}
S_j \frac{\partial}{\partial s_j}f(\mathbf{S}) 
&=& \sum_{i=1}^m \xi_{i,j} \frac{\partial}{\partial s_j}f(\mathbf{S}) \nonumber \\
&=& \sum_{i=1}^m \xi_{i,j}\frac{\partial}{\partial s_j}f(S_1,\ldots, S_{j*i} , \ldots, S_p) \nonumber \\ 
&\text{ \ \ \ }& + \sum_{i=1}^m \xi_{i,j}^2 \int_0^1 \frac{\partial ^2}{\partial s_j^2} f(S_1,\ldots, S_{j*i}+u\xi_{i,j} , \ldots ,S_p) du .
\ena

Substituting (\ref{2nd}) and (\ref{3rd}) into (\ref{1st}) and letting $\mathbf{S}^{j*i} =: (S_1,\ldots,S_{j-1},S_{j*i},S_{j+1},\ldots,S_p)$, we obtain
\bea \label{4th}
\e\left[h(\Sigma^{-1/2}\mathbf{S}) - N h \right] 
&=& \e\sum_{j=1}^p \sum_{i=1}^m \xi_{i,j}^2 \int_0^1 \left( \frac{\partial ^2}{\partial s_j^2}f(\mathbf{S}) -\frac{\partial ^2}{\partial s_j^2} f(S_1,\ldots, S_{j*i}+u\xi_{i,j} , \ldots ,S_p) \right) du \nonumber \\
&\text{ \ \ }&  + \e\sum_{j=1}^p \sum_{i=1}^m (\sigma_{i,j}^2 - \xi_{i,j}^2) \frac{\partial ^2}{\partial s_j^2}f(\mathbf{S}) - \e \sum_{j=1}^p\sum_{i=1}^m \xi_{i,j}\frac{\partial}{\partial s_j}f(\mathbf{S}^{j*i}) \nonumber \\
&\text{ \ \ }& + \e\sum_{j=1}^p\sum_{i,k \in [m],i \neq k} \sigma_{i,j;k,j}\frac{\partial ^2}{\partial s_j^2}f(\mathbf{S}) + \e\sum_{j,l \in [p],j \ne l} \Sigma_{j,l}\frac{\partial ^2}{\partial s_j\partial s_l}f(\mathbf{S}).
\ena

Now we handle these five terms in (\ref{4th}) separately. For the first term, using \eqref{diffbound} we have
\bea \label{5th}
\Bigg|\e\sum_{j=1}^p \sum_{i=1}^m &\xi_{i,j}^2& \int_0^1 \left( \frac{\partial ^2}{\partial s_j^2}f(\mathbf{S}) -\frac{\partial ^2}{\partial s_j^2} f(S_1,\ldots, S_{j*i}+u\xi_{i,j} , \ldots ,S_p) \right) du\Bigg| \nonumber \\
&=& \left|\e\sum_{j=1}^p  \sum_{i=1}^m \xi_{i,j}^2 \int_0^1 \int_{u\xi_{i,j}}^{\xi_{i,j}}\frac{\partial ^3}{\partial s_j^3} f(S_1, \ldots , S_{j*i}+t, \ldots,S_p)dtdu\right| \nonumber \\
&\le& \frac{p^3}{3}|\Sigma^{-1/2}|_{\infty}^3\e\sum_{j=1}^p  \sum_{i=1}^m \xi_{i,j}^2 \int_0^1 \int_{u|\xi_{i,j}|}^{|\xi_{i,j}|}dtdu \nonumber \\
&=& \frac{p^3}{6}|\Sigma^{-1/2}|_{\infty}^3\sum_{j=1}^p  \sum_{i=1}^m \e|\xi_{i,j}|^3 \nn \\
&\le& \frac{p^3}{6}|\Sigma^{-1/2}|_{\infty}^3B\sum_{j=1}^p  \sum_{i=1}^m \e\xi_{i,j}^2 \nn \\
&=& \frac{p^3}{6}|\Sigma^{-1/2}|_{\infty}^3 B \sum_{j=1}^p\left(\Sigma_{j,j}- \sum_{i,k \in [m],i \neq k} \sigma_{i,j;k,j}\right),
\ena
where we have used the almost sure bound on the variables $\xi_{i,j}$, and that their sum $S_j$ over $i$ from $1$ to $m$ has mean zero in the last two inequalities, respectively.

For the second term in (\ref{4th}), we first bound it by the three terms denoted by $I_1,I_2,I_3$, respectively, 
\beas 
&&\left|\e\sum_{j=1}^p \sum_{i=1}^m (\sigma_{i,j}^2 - \xi_{i,j}^2) \frac{\partial ^2}{\partial s_j^2}f(\mathbf{S}) \right| \nn\\
&&\hspace{30pt}\le  \left|\e\sum_{j=1}^p \sum_{i=1}^m (\sigma_{i,j}^2 - \xi_{i,j}^2) \frac{\partial ^2}{\partial s_j^2}f(\mathbf{S}^{j*i}) \right|
+\left|\e\sum_{j=1}^p \sum_{i=1}^m (\sigma_{i,j}^2 - \xi_{i,j}^2) \left(\frac{\partial ^2}{\partial s_j^2}f(\mathbf{S}^{j*i})-\frac{\partial ^2}{\partial s_j^2}f(\mathbf{S})\right) \right| \nn\\
&&\hspace{30pt}\le \left|\sum_{j=1}^p \sum_{i=1}^m  \Cov\left(\frac{\partial ^2}{\partial s_j^2}f(\mathbf{S}^{j*i}),\xi_{i,j}^2\right)\right| 
+\e\sum_{j=1}^p \sum_{i=1}^m \sigma_{i,j}^2  \left|\left(\frac{\partial ^2}{\partial s_j^2}f(\mathbf{S}^{j*i})-\frac{\partial ^2}{\partial s_j^2}f(\mathbf{S})\right) \right| \nn\\
&&\hspace{40pt}+\e\sum_{j=1}^p \sum_{i=1}^m \xi_{i,j}^2  \left|\left(\frac{\partial ^2}{\partial s_j^2}f(\mathbf{S}^{j*i})-\frac{\partial ^2}{\partial s_j^2}f(\mathbf{S})\right) \right| := I_1+I_2+I_3 .
\enas
Then we write $I_1$ as 
\beas
I_1=\left|\sum_{j=1}^p \sum_{i=1}^m  \Cov\left(\frac{\partial ^2}{\partial s_j^2}f(\mathbf{S}^{j*i}),g(\xi_{i,j})\right)\right|,
\enas
where 
\beas
g(x)=\left\{  
\begin{array}{cc}
x^2 & |x| \le B\\
B^2 & |x|>B,
\end{array} \right.
\enas
As ${\bf S}^{j*i}$ and $\xi_{i,j}$ are increasing functions defined on disjoint subsets of $\bsxi$, by Lemma \ref{lemma:inc}, $({\bf S}^{j*i},\xi_{i,j})$ are negatively associated for all $i,j$.  Applying Lemma \ref{lemma:covbound} and using the bound \eqref{diffbound}, we obtain
\bea \label{i1}
I_1 
&\le&  \left|\sum_{j=1}^p \sum_{i=1}^m \left(\sum_{l \in [p],l \ne j}  \left|\frac{\partial^3}{\partial s_l\partial s_j^2} f \right|_{\infty}\left|\frac{\partial g}{\partial x}\right|_{\infty} \Cov\left(S_l,\xi_{i,j}\right) + \left|\frac{\partial^3}{\partial s_j^3} f \right|_{\infty}\left|\frac{\partial g}{\partial x}\right|_{\infty} \Cov\left(S_{j*i},\xi_{i,j}\right)\right) \right| \nonumber \\
&\le& \frac{2}{3}p^3|\Sigma^{-1/2}|_{\infty}^3B\left|\sum_{j=1}^p \sum_{i=1}^m \left(\sum_{l \in [p],l \ne j}  \Cov\big(S_l,\xi_{i,j}\big) + \Cov\left(S_{j*i},\xi_{i,j}\right)\right)\right| \nonumber \\
&=&  \frac{2}{3}p^3|\Sigma^{-1/2}|_{\infty}^3B\left|\sum_{j,l \in [p], j\ne l} \sum_{i,k=1}^m   \sigma_{i,j;k,l} + \sum_{j=1}^p \sum_{i,k\in[m],i\ne k} \sigma_{i,j;k,j}\right| \nonumber \\
&=&  -\frac{2}{3}p^3|\Sigma^{-1/2}|_{\infty}^3B\left( \sum_{j,l \in [p], j\ne l} \Sigma_{j,l}+ \sum_{j=1}^p \sum_{i,k\in[m],i\ne k} \sigma_{i,j;k,j}\right).
\ena
Again using \eqref{diffbound}, we have
\bea
I_2 &\le& \sum_{j=1}^p \sum_{i=1}^m \sigma_{i,j}^2 \left|\frac{\partial ^3}{\partial s_j^3}f\right|_{\infty} \e\left|S_{j*i}-S_j \right| 
   = \sum_{j=1}^p \sum_{i=1}^m \sigma_{i,j}^2 \left|\frac{\partial ^3}{\partial s_j^3}f\right|_{\infty} \e\left|\xi_{i,j} \right| \nn\\
	&\le& \frac{p^3}{3}|\Sigma^{-1/2}|_{\infty}^3B \sum_{j=1}^p \left(\Sigma_{j,j} -\sum_{i,k\in[m],i\ne k} \sigma_{i,j;k,j} \right),
\ena 
and
\bea
I_3 &\le&  \sum_{j=1}^p \sum_{i=1}^m \left|\frac{\partial ^3}{\partial s_j^3}f\right|_{\infty} \e \xi_{i,j}^2\left|S_{j*i}-S_j \right|  
\le \frac{p^3}{3}|\Sigma^{-1/2}|_{\infty}^3B \sum_{j=1}^p \sum_{i=1}^m \e\xi_{i,j}^2  \nn \\
&=&\frac{p^3}{3}|\Sigma^{-1/2}|_{\infty}^3B \sum_{j=1}^p \left(\Sigma_{j,j} -\sum_{i,k\in[m],i\ne k} \sigma_{i,j;k,j} \right).
\ena

For the third term in (\ref{4th}), again applying Lemma \ref{lemma:covbound} and arguing as for $I_1$ in \eqref{i1}, we have
\bea 
&&\left|\e\sum_{j=1}^p \sum_{i=1}^m  \xi_{i,j}\frac{\partial}{\partial s_j}f(\mathbf{S}^{j*i})\right| 
= \left|\sum_{j=1}^p \sum_{i=1}^m \Cov\left(\xi_{i,j},\frac{\partial}{\partial s_j}f(\mathbf{S}^{j*i}) \right)\right| \nonumber \\
&&\text{ \ \ \ \ \ \ \  } \le \left|\sum_{j=1}^p \sum_{i=1}^m \left(\sum_{l \in [p]\setminus \{j\}}  \left|\frac{\partial^2}{\partial s_l\partial s_j}f\right|_{\infty}  \Cov\left(\xi_{i,j},S_l\right) + \left|\frac{\partial^2}{\partial s_j^2}f\right|_{\infty}\Cov\left(\xi_{i,j},S_{j*i} \right) \right)\right| \nonumber \\
&&\text{ \ \ \ \ \ \ \  } \le -\frac{p^2}{2}|\Sigma^{-1/2}|_{\infty}^2\left( \sum_{j,l \in [p],j \ne l}\sum_{i=1}^m\Cov\left(\xi_{i,j},S_l\right)  + \sum_{j=1}^p \sum_{i=1}^m \Cov \left(\xi_{i,j},S_{j*i} \right)\right) \nonumber \\
&&\text{ \ \ \ \ \ \ \  } = -\frac{p^2}{2}|\Sigma^{-1/2}|_{\infty}^2 \left(\sum_{j,l \in [p],j \ne l}\sum_{i,k=1}^m \sigma_{i,j;k,l} +  \sum_{j=1}^p\sum_{i,k \in [m],i \neq k} \sigma_{i,j;k,j} \right) \nonumber \\
&&\text{ \ \ \ \ \ \ \  } = -\frac{p^2}{2}|\Sigma^{-1/2}|_{\infty}^2 \left(\sum_{j,l \in [p],j \ne l}\Sigma_{j,l} +  \sum_{j=1}^p\sum_{i,k \in [m],i \neq k} \sigma_{i,j;k,j} \right).
\ena

For the fourth and the fifth terms in (\ref{4th}), again using \eqref{diffbound} we obtain
\bea
\left|\sum_{j=1}^p\sum_{i,k \in [m],i \neq k} \sigma_{i,j;k,j}\frac{\partial ^2}{\partial s_j^2}f(\mathbf{S})\right| \le -\frac{p^2}{2}|\Sigma^{-1/2}|_{\infty}^2\sum_{j=1}^p\sum_{i,k \in [m],i \neq k} \sigma_{i,j;k,j}
\ena
and
\bea \label{9th}
\left|\e\sum_{j,l \in [p],j \ne l} \Sigma_{j,l}\frac{\partial ^2}{\partial s_j\partial s_l}f(\mathbf{S})\right| \le  -\frac{p^2}{2}|\Sigma^{-1/2}|_{\infty}^2\sum_{j,l \in [p],j \ne l} \Sigma_{j,l}.
\ena
Summing the bounds (\ref{5th})-(\ref{9th}) we find that $\left|\e\left[h(\Sigma^{-1/2}\mathbf{S}) - N h \right]\right|$ is bounded by the right hand side of \eqref{steinbound:fkg2}. Taking supremum over $h \in \mathcal{H}_{3,\infty,p}$ and using the definition \eqref{smoothdef} of $d_{\mathcal{H}_{m,\infty,p}}$, yields the claim.

\bbox

\section{Applications} \label{app}

In this section, we follow the same structure as in Section 2.1 of \cite{GW16} with positive association replaced by negative association. In particular, we apply our main theorems in the first section to second order stationary negatively associated random fields with exponential covariance decay.

First we introduce the definitions and notations used in \cite{GW16} that will also be used here. Let $\{ X_{{\bf j}}: {\bf j}\in \mathbb{Z}^d \}$ be a negatively associated random field on the $d$-dimensional integer lattice $\mathbb{Z}^d$ and assume that the field is second order stationary. We recall that a random field $\{X_{{\bf j}}: {\bf j} \in \mathbb{Z}^d\}$ is called {\em second order stationary} when $\e X_{{\bf j}}^2<\infty$ for all ${\bf j} \in \mathbb{Z}^d$ and the covariance $\Cov \big(X_{{\bf i}},X_{{\bf j}} \big) = R({\bf j}-{\bf i})$ for all ${\bf i},{\bf j}\in \mathbb{R}^d$, with $R(\cdot)$ given by
\bea \label{covdef}
R({\bf k})=\Cov \big(X_{{\bf 0}},X_{{\bf k}} \big).
\ena

We let ${\bf 1} \in \mathbb{Z}^d$ denotes the vector with all components $1$, and write inequalities such as ${\bf a} < {\bf b}$ for vectors ${\bf a} , {\bf b} \in \mathbb{R}^d$  when they hold componentwise.
For ${\bf k} \in \mathbb{Z}^d, n \in \mathbb{N}_1$, define the `block sum' variables, over a block with side length $n$, by 
\bea \label{blockdef}
S_{{\bf k}}^n =  \sum_{{\bf j} \in B_{{\bf k}}^n} X_{{\bf j}} \qmq{where}
B_{{\bf k}}^n = \left\{ {\bf j} \in \mathbb{Z}^d: {\bf k} \le {\bf j} < {\bf k}+n{\bf 1} \right\}.
\ena
Note that $B_{\bf k}^n=B_{\bf 0}^n+{\bf k}$.

For $R(\cdot)$ given in \eqref{covdef}, we have
\bea \label{andef}
\Var \big(S_{{\bf k}}^n \big)  = \sum_{{\bf i},{\bf j} \in B_{{\bf k}}^n} \Cov \big(X_{{\bf i}},X_{{\bf j}} \big) = n^d A_n \qmq{where} A_n = \frac{1}{n^d} \sum_{{\bf i},{\bf j} \in B_{{\bf 1}}^n} R({\bf i}-{\bf j}).
\ena 
Since the field is negatively associated, $R({\bf k}) \le 0$ for all ${\bf k} \ne {\bf 0},$ which implies that $0\le A_n \le R({\bf 0})$ for all $n \in \mathbb{N}_1$. For simplicity, in this work, we assume that $\inf_n A_n > 0$ which implies that $A_n$ is of constant order. With this assumption, we may include $A_n$ in our bounds without affecting the rate of convergence.

With $S_{\bf k}^n$ defined in \eqref{blockdef}, we consider the standardized variables
\bea \label{wn2}
W_{{\bf k}}^n = \frac{S_{{\bf k}}^n - \e S_{{\bf k}}^n}{\sqrt{n^dA_n}}, \qmq{${\bf k} \in \mathbb{Z}^d, n \in \mathbb{N}_1$,}
\ena
that have mean zero and variance 1. The following theorem provides a bound of order $n^{-d/(2d+2)}$ with an explicit constant
on the $L^1$ distance between the distribution of $W_{{\bf k}}^n$ and the normal under the assumption that the covariance function $R(\cdot)$ decays at exponential rate in the $L^1$ norm in $\mathbb{R}^d$. Since all norms in $\mathbb{R}^d$ are equivalent, we use the $L^1$ norm that makes our calculation simplest. 

\begin{theorem} \label{1stappthm} 
Let $d \in \mathbb{N}_1$ and  $\{ X_{{\bf j}}: {\bf j}\in \mathbb{Z}^d \}$ be a negatively associated second order stationary random field with covariance function $R({\bf k})={\rm Cov}(X_{{\bf j}},X_{{\bf j}+{\bf k}})$ for all ${\bf j},{\bf k} \in {\mathbb Z}^d$, and 
suppose that for some $K > 0$, $|X_{{\bf j}}| \le K$ a.s.\ for all ${\bf j} \in \mathbb{Z}^d$. Assume that there exist $\lambda >0$ and $\kappa_0 >0$ such that 
\bea \label{Ck.exp.bound2}
-R({\bf k}) \le \kappa_0 e^{-\lambda |{\bf k}|_1} \qmq{for all ${\bf k} \in \mathbb{Z}^d/\{\mathbf{0}\}$ }
\ena
and $\inf_n A_n > 0$ where $A_n$ is given in \eqref{andef}. Let
\bea \label{muandnu}
\mu_{\lambda} = \frac{e^{\lambda}}{\left(e^{\lambda}-1\right)^2} \text{, \ } \nu_{\lambda} = \frac{e^{2\lambda}}{\left(e^{\lambda}-1\right)^2} \qmq{and} \gamma_{\lambda,d}= (4\mu_\lambda+2\nu_\lambda)^d-(2\nu_\lambda)^d
\ena	
and	
\bea \label{cdef}
C_{\lambda,\kappa_0,d}=\frac{10Kd\sqrt{A_n}}{5.2\kappa_0\gamma_{\lambda,d}}
\ena 
Then, for any ${\bf k} \in \mathbb{Z}^d$, with $W_{{\bf k}}^n$ as given in \eqref{wn2} and $Z$ a standard normal random variable, 
\bea \label{1stappbound}
d_1\big({\cal L}(W_{{\bf k}}^n),{\cal L}(Z)\big) \le \frac{\kappa_1}{n^{d/(2d+2)}} 
\text{ \ \  for all \ \  } n \ge \max \left\{ C_{\lambda,\kappa_0,d}^{2/d},C_{\lambda,\kappa_0,d}^{-2/(d+2)} \right\}, 
\ena
where
\beas
\kappa_1 =  \left(\frac{10K (5.2\kappa_0 \gamma_{\lambda,d})^d }{A_n^{d+1/2}}\right)^{1/(d+1)} \left(\frac{1
}{d^{\frac{d}{d+1}}}+ 2d^{\frac{1}{d+1}}\right). 
\enas
\end{theorem}

The bound in \eqref{1stappbound} is of order $n^{-d/(2d+2)}$ recalling that $A_n$ is bounded away from zero and infinity and hence does not contribute the rate.
We also extend Theorem \ref{1stappthm} to the multidimensional case. For any $p \in \mathbb{N}_1$ and indices ${\bf k}_1,\ldots,{\bf k}_p$ in $\mathbb{Z}^d$  such that $B_{{\bf k}_i}^n, i \in [p]$ are disjoint, Theorem \ref{1stappthm2} provides a bound in the metric $d_{\mathcal{H}_{3,\infty,p}}$ to the multivariate normal for ${\bf S}^n=(S_{{\bf k}_1}^n,\ldots,S_{{\bf k}_p}^n)$ under exponential decay of the covariance function. We note the difference between the assumption \eqref{eq:kqks.separated.alphan} below and the one in (30) of \cite{GW16} that we require $B_{{\bf k}_i}^n, i \in [p]$ here to be disjoint, otherwise, there exists ${\bf j} \in \mathbb{Z}^d$ that belongs to both $B_{{\bf k}_i}^n$ and $B_{{\bf k}_j}^n$ for some $i,j$ and $(X_{\bf j},X_{\bf j})$ is positively associated to which Theorem \ref{stein:fkg2} is not applicable. The same issue does not arise in the positively associated case as a pair of the same variable has positive covariance.  In the following result and its proof, constants will not be tracked with precision, but will be indexed by the set of variables on which it depends.

\begin{theorem} \label{1stappthm2}
For $d \in \mathbb{N}_1$, let $\{ X_{{\bf j}}: {\bf j}\in \mathbb{Z}^d \}$ be a negatively associated second order stationary random field with covariance function $R({\bf k})={\rm Cov}(X_{{\bf j}},X_{{\bf j}+{\bf k}})$ for all ${\bf j},{\bf k} \in {\mathbb Z}^d$, and suppose that there exist constants $K>0,\kappa_0>0$ and $\lambda > 0$ such that $|X_{{\bf j}}| \le K$ a.s.\ for all ${\bf j} \in \mathbb{Z}^d$,
\beas
-R({\bf k}) \le \kappa_0e^{-\lambda |{\bf k}|_1} \text{ \ for all \ } {\bf k}\in \mathbb{Z}^d/\{ \mathbf{0}\},
\enas
and $\inf_n A_n > 0$ where $A_n$ is given in \eqref{andef}.
For $p\in \mathbb{N}_1$ let ${\bf k}_1,\ldots,{\bf k}_p \in \mathbb{Z}^d$ be 
such that
\bea \label{eq:kqks.separated.alphan}
\min_{q,s\in[p], q \ne s}\left|{\bf k}_q-{\bf k}_s\right|_{\infty} \ge n.
\ena
Let
$\mathbf{S}^n = (S_{{\bf k}_1}^n ,\ldots,S_{{\bf k}_p}^n )$, where $S_{{\bf k}}^n$ is defined as in \eqref{blockdef} and $\Sigma$ be the covariance matrix of $\mathbf{S}^n$. Then, for $n > (p-1)\kappa_0\nu_{\lambda}^de^{-\lambda}A_n^{-1}$ with $\nu_\lambda$ as in \eqref{muandnu}, $\Sigma$ is invertible and
\bea \label{eq:infinity.norm.bound.inverse.lem:Gershgorin}
|\Sigma^{-1}|_{\infty} \le \frac{1}{n^{d-1}(nA_n - (p-1)\kappa_0 \nu_{\lambda}^d e^{-\lambda})}.
\ena
Furthermore, with
\beas 
\psi_n=n^{d/2}|\Sigma^{-1/2}|_\infty \text{ \ for \ } n > (p-1)\kappa_0\nu_{\lambda}^de^{-\lambda}A_n^{-1} \text{ \ \ and \ \ } B_{n,d}= d\psi_n A_n
\enas 
there exists a constant $C_{\lambda,\kappa_0,d,p,K}$ such that, for  
\bea \label{eq:thm2rangen}
n > \max \left\{B_{n,d}^{2/d},B_{n,d}^{-2/(d+2)},(p-1)\kappa_0\nu_{\lambda}^de^{-\lambda}A_n^{-1}\right\},
\ena
\bea\label{2ndappbound}
&&d_{\mathcal{H}_{3,\infty,p}}\big(\mathcal{L}(\Sigma^{-1/2}(\mathbf{S}^n-\e \mathbf{S}^n)),\mathcal{L}(\mathbf{Z})\big) \nn\\ 
&&\hspace{60pt}\le C_{\lambda,\kappa_0,d,p,K} \Bigg(\frac{\psi_n^{(2d+4)/(d+1)}}{A_n^{(d-1)/(d+1)}n^{d/(d+1)}} 
+ \frac{\psi_n^{(2d+3)/(d+1)}}{A_n^{d/(d+1)}n^{(3d+2)/(2d+2)}}  \nn\\ 
&&\hspace{150pt}+ \frac{ \psi_{n}^2}{n} +\frac{A_n^{1/(d+1)}\psi_n^{(2d+3)/(d+1)}}{n^{d/(2d+2)}}  \Bigg),
\ena
where $\mathbf{Z}$ is a standard normal random vector in $\mathbb{R}^p$.
\end{theorem}

Since $A_n$ is of constant order, $|\Sigma^{-1/2}|_\infty$ is of order at most $n^{-d/2}$ by \eqref{eq:infinity.norm.bound.inverse.lem:Gershgorin}. This implies that $\psi_n$ is of at most constant order and thus so is $B_{n,d}$. Therefore the last term on the right hand side of \eqref{2ndappbound} is the only one that contributes the rate of convergence of order $n^{-d/(2d+2)}$ as the other terms converge to zero at much faster rates. We note that the bounds in Theorems \ref{1stappthm} and \ref{1stappthm2} have the same order as the ones in Theorems 2.1 and 2.2 of \cite{GW16}, respectively. However, comparing to the results of \cite{GW16}, the constant of the bound \eqref{1stappbound} is bigger and the bound \eqref{2ndappbound} has the extra terms that do not contribute the rate.

To prove Theorems \ref{1stappthm} and \ref{1stappthm2}, we use the same technique as in \cite{GW16} decomposing the sum $S_{\bf k}^n$ over the block $B_{\bf k}^n$
into sums over smaller, disjoint blocks whose side lengths are at most some integer $l$. That is, 
for $1 \le l \le n$, we uniquely write $n=(m-1)l+r$ with $m \ge 1$ and $1 \le r \le l$ and correspondingly decompose $B_{{\bf k}}^n$ into $m^d$ disjoint blocks $D_{{\bf i},{\bf k}}^l, {\bf i} \in [m]^d$, where there are $(m-1)^d$ `main' blocks having all sides of length $l$, and $m^d-(m-1)^d$ remainder blocks having all sides of length $r$ or $l$, with at least one side of length $r$.

To be more precise, for ${\bf k} \in \mathbb{Z}^d$ and ${\bf i} \in [m]^d$ set $D_{{\bf i},{\bf k}}^l=D_{\bf i}^l+{\bf k} - {\bf 1}$ where
\beas
D_{\bf i}^l&=&\big\{{\bf j} \in \mathbb{Z}^d: (i_s-1)l+1 \le j_s \le i_sl \text{  for  } i_s \not =m, \\
                                                   &&\hspace{50pt}(m-1)l+1 \le j_s \le (m-1)l+r \text{  for  } i_s =m \big\}.
\enas
It is easy to see that for ${\bf i} \in [m-1]^d$, the vectors indexing the `main blocks', we have 
\bea \label{Di.block}
D_{\bf i}^l = B_{({\bf i}-{\bf 1})l+{\bf 1}}^l  \text{ \ \ for \ \ } {\bf i} \in [m-1]^d,
\ena
and if $r=l$ then $D_{{\bf i}}^l$ is given by \eqref{Di.block} for all $i \in [m]^d$. Furthermore, it is straightforward to verify that the elements of the collection $\{D_{{\bf i},{\bf k}}^l, {\bf i} \in [m]^d\}$ is the partition of $B_{\bf k}^n$.

Letting
\bea \label{xidef}
\xi_{{\bf i},{\bf k}}^{l} = \sum_{{\bf t} \in D_{{\bf i},{\bf k}}^l}\left(X_{{\bf t}} - \e X_{{\bf t}}\right) \text{ \ \ for \ \  } {\bf i} \in [m]^d, \quad \mbox{and} \quad W_{{\bf k}}^n = \sum_{{\bf i} \in [m]^d} \frac{\xi_{
		{\bf i},{\bf k}
		}^{l}}{\sqrt{n^d A_n}} ,
\ena 
we see that $\xi_{{\bf i},{\bf k}}^n$ has mean zero, and
$W_{{\bf k}}^n$ as in \eqref{xidef} agrees with its representation as given in \eqref{wn2}, and has mean zero and variance one. For simplicity we will drop the index ${\bf k}$ in $\xi_{{\bf i},{\bf k}}$ when ${\bf k}={\bf 1}$, as we do also for $D_{{\bf i},{\bf k}}$, and also suppress $n$ in $\xi_{{\bf i},{\bf k}}^n$. 

As the elements of $\{ \xi_{{\bf i},{\bf k}}: {\bf i} \in [m]^d \}$ are increasing functions of disjoint subsets of $\{ X_{{\bf j}}: {\bf j} \in \mathbb{Z}^d \}$, they are negatively associated by Lemma \ref{lemma:inc}. 
We prove Theorems \ref{1stappthm}  and \ref{1stappthm2} with the help of the following three lemmas. The first, Lemma \ref{lemma:1stapp} bounds the sum of the covariances between $\xi_{{\bf i},{\bf k}}^{l}$ and $\xi_{{\bf j},{\bf k}}^{l}$, defined in \eqref{xidef}, over ${\bf i},{\bf j} \in [m]^d$. Next, we state Lemma \ref{sumexpo}, proved in \cite{GW16}, which is used in the proof of Lemma \ref{lemma:1stapp2} that bounds the covariance between two non-overlapped block sums of size $n^d$. 

\begin{lemma} \label{lemma:1stapp}
Let $\{ X_{{\bf j}}: {\bf j}\in \mathbb{Z}^d \}$ be a second order, negatively associated stationary random field with covariance function $R({\bf k})={\rm Cov}(X_{{\bf j}},X_{{\bf j}+{\bf k}})$ for all ${\bf j} \in {\mathbb Z}^d$ and ${\bf k} \in {\mathbb Z}^d/\{\mathbf{0}\}$ where $R(\cdot)$ satisfies the exponential decay condition \eqref{Ck.exp.bound2}. For $n \ge 2$ and $1 \le l \le n$, let $n=(m-1)l+r$ for integers $m \in \mathbb{N}_1$ and $1 \le r \le l$. Then for ${\bf i} \in [m]$ and ${\bf k} \in \mathbb{Z}^d$, with $\xi_{{\bf i},{\bf k}}$ given by \eqref{xidef} we have 
\beas
\sum_{{\bf i},{\bf j} \in [m]^d,{\bf i} \neq {\bf j}} -\e \left[\xi_{{\bf i},{\bf k}}^{l} \xi_{{\bf j},{\bf k}}^{l}\right] \le \frac{\kappa_0\gamma_{\lambda,d}n^d}{l},
\enas
where $\kappa_0$ is given in \eqref{Ck.exp.bound2}, and $\gamma_{\lambda,d}$ in \eqref{muandnu}. 
\end{lemma}

\proof We note that the difference between \eqref{Ck.exp.bound2} of this work and (22) of \cite{GW16} is only the sign on the left hand side of the inequality. Thus the proof follows immediately from the proof of Lemma 2.4 of \cite{GW16} with $\e \left[\xi_{{\bf i},{\bf k}}^{l} \xi_{{\bf j},{\bf k}}^{l}\right]$ and $R({\bf k})$ replaced by $-\e \left[\xi_{{\bf i},{\bf k}}^{l} \xi_{{\bf j},{\bf k}}^{l}\right]$ and $-R({\bf k})$, respectively.

\bbox

\begin{lemma}[\cite{GW16}] \label{sumexpo} 
For all $n \in \mathbb{N}_2$ and $\lambda > 0$,
\beas
\sum_{a=-n+1}^{n-1} (n-|a|)e^{-\lambda |q+a|}  \text{ \ is decreasing as a function of $|q| \in \mathbb{N}_0$. }
\enas
\end{lemma}

In the following we will use the identities
\bea \label{Sum.krk4}
\sum_{k=1}^{n-1} (n-k)w^{k}= \frac{w\left((n-1)-nw+w^n\right)}{(w-1)^2} \text{ \ \ for \ \ } w \ne 1,
\ena
and
\bea \label{Sum.krk2}
n+ 2\sum_{b=1}^{n-1} (n-b)u^{b}= \frac{(1-u^2)n-2u+2u^{n+1}}{(u-1)^2} \text{ \ \ for \ \ } u \ne 1.
\ena

\begin{lemma} \label{lemma:1stapp2}
Let $\{ X_{{\bf j}}: {\bf j}\in \mathbb{Z}^d \}$ be a second order stationary random field with covariance function $R({\bf k})={\rm Cov}(X_{{\bf j}},X_{{\bf j}+{\bf k}})$ for all ${\bf j} \in {\mathbb Z}^d$ and ${\bf k} \in {\mathbb Z}^d/\{\mathbf{0}\}$ where $R(\cdot)$ satisfies \eqref{Ck.exp.bound2}. Let $n \in \mathbb{N}_2$, ${\bf k}_1$ and ${\bf k}_2$ be vectors in $\mathbb{Z}^d$ such that
\beas
\left|{\bf k}_1-{\bf k}_2\right|_{\infty} \ge n,
\enas
then with $\lambda$ and $\kappa_0$ as in \eqref{Ck.exp.bound2}, and $\nu_\lambda$ as in \eqref{muandnu},
\beas
-\Cov \left(S_{{\bf k}_1}^n,S_{{\bf k}_2}^n \right) \le \kappa_0 \nu_\lambda^d e^{-\lambda} n^{d-1}.
\enas
\end{lemma}

\proof
Using the definition of $S_{\bf k}^n$ and that $X_{{\bf j}}, {\bf j}\in \mathbb{Z}^d$ are second order stationary, we have
\bea
-\Cov \left(S_{{\bf k}_1}^n,S_{{\bf k}_2}^n \right)&=&-\sum_{\substack{p_1,\ldots,p_d =0 \\  q_1,\ldots,q_d=0}}^{n-1} R\left(\begin{bmatrix} (p_1+k_1^2)-(q_1+k_1^1) \\ \vdots \\ (p_d+k_d^2)-(q_d+k_d^1) \end{bmatrix}\right) \nonumber \\
&\le& \kappa_0 \sum_{a_1,\ldots,a_d=-n+1}^{n-1} (n-|a_1|)\ldots(n-|a_d|) \exp  \left(-\lambda\left|\begin{bmatrix} a_1+(k_1^2-k_1^1)\\ \vdots \\ a_d+(k_d^2-k_d^1) \end{bmatrix}\right|_1\right) \nonumber \\
&=& \kappa_0 \prod_{i=1}^d \sum_{a_i=-n+1}^{n-1} (n-|a_i|)e^{-\lambda |(k_i^2-k_i^1)+a_i|}.\label{eq:prod.term.one.exception}
\ena
Applying Lemma \ref{sumexpo}, we have 
\bea \label{eq:dec.in.k1-k2}
\sum_{a_i=-n+1}^{n-1} (n-|a_i|)e^{-\lambda |(k_i^2-k_i^1)+a_i|} \qmq{is a decreasing function of $|k_i^1-k_i^2|$.}
\ena
Hence the $i^{th}$ sum appearing in the product
\eqref{eq:prod.term.one.exception}	
is maximized by its value when $k_i^1=k_i^2$. As $\left|{\bf k}_1-{\bf k}_2\right|_{\infty} \ge n$, there must exist at least one $i$ for which  $|k_i^2-k_i^1| \ge n$, and whose corresponding sum is bounded by its value when $|k_i^2-k_i^1|$ is exactly $n$, using \eqref{eq:dec.in.k1-k2}. The product of these sums, by \eqref{eq:dec.in.k1-k2} again, is maximized when there is just a single coordinate achieving $n$ as its absolute difference, and where this difference in all other terms achieve equality to zero.  Therefore, by symmetry \eqref{eq:prod.term.one.exception} is bounded by the case where $k_i^1=k_i^2$ for $i \in [d-1]$ and $k_d^2-k_d^1=n$ and thus
\bea 
-\Cov \left(S_{{\bf k}_1}^n,S_{{\bf k}_2}^n \right) 
&\le& \kappa_0\prod_{i=1}^{d-1} \sum_{a_i=-n+1}^{n-1} (n-|a_i|) e^{-\lambda|a_i|} \sum_{a_d=-n+1}^{n-1} (n-|a_d|)e^{-\lambda|a_d+n|} \nn\\
&\le& \kappa_0\left(n\nu_{\lambda} \right)^{d-1} \sum_{a_d=-n+1}^{n-1} (n-|a_d|)e^{-\lambda |a_d+n|},\label{sumlemma2.5}
\ena
where we have applied \eqref{Sum.krk2} in the final inequality and $\nu_{\lambda}$ is given in \eqref{muandnu}.

Now considering the sum in \eqref{sumlemma2.5}, we obtain
\beas
\sum_{a=-n+1}^{n-1} (n-|a|)e^{-\lambda |a+n|} &=& \sum_{a=-n+1}^{n-1} (n-|a|)e^{-\lambda (a+n)} \\
&=& \sum_{a=1}^{n-1}(n-a)e^{-\lambda (-a+n)} + \sum_{a=0}^{n-1}(n-a)e^{-\lambda (a+n)}.
\enas 
For each sum, making a change of variable and applying \eqref{Sum.krk4}, we obtain 
\beas
\sum_{a=1}^{n-1} (n-a)e^{-\lambda (-a+n)} &=& e^{-\lambda n}\sum_{a=1}^{n-1} (n-a)e^{\lambda a} \\
&=& \frac{e^{-\lambda n}e^{\lambda}\left(n-1-ne^{\lambda}+e^{\lambda n}\right)}{\left(e^{\lambda}-1\right)^2}\\
&=& \frac{e^{\lambda}+n\left(1-e^{\lambda}\right)e^{\lambda(1-n)}-e^{\lambda(1-n)}}{\left(e^{\lambda}-1\right)^2}
\enas
and
\beas
\sum_{a=0}^{n-1} (n-a)e^{-\lambda (a+n)} &=&ne^{-\lambda n}+e^{-\lambda n}\sum_{a=1}^{n-1} (n-a)e^{-\lambda a} \\
&=& ne^{-\lambda n}+\frac{e^{-\lambda n}e^{-\lambda}\left(n-1-ne^{-\lambda}+e^{-\lambda n}\right)}{\left(e^{-\lambda}-1\right)^2} \\
&=& ne^{-\lambda n}+\frac{e^{-\lambda n}e^{\lambda}\left(n-1-ne^{-\lambda}+e^{-\lambda n}\right)}{\left(e^{\lambda}-1\right)^2} \\
&=&  \frac{e^{\lambda(1-2n)}-e^{\lambda(1-n)}+n(e^{\lambda}-1)e^{\lambda(1-n)}}{\left(e^{\lambda}-1\right)^2}.
\enas

Summing these two terms yields
\beas
\sum_{a=-n+1}^{n-1} (n-|a|)e^{-\lambda |a+n|}
=\frac{e^{\lambda(1-2n)}+e^{\lambda}-2e^{\lambda(1-n)}}{\left(e^{\lambda }-1\right)^2}
\le \frac{e^{ \lambda }}{\left(e^{\lambda }-1\right)^2} =  \nu_{\lambda} e^{-\lambda}.
\enas
Plugging the last bound in \eqref{sumlemma2.5} yields the claim.

\bbox

Now we have all ingredients to prove Theorems \ref{1stappthm} and \ref{1stappthm2}. In the following, we use the same technique as in (44) of \cite{GW16}, that is, for any positive real numbers $a$ and $b$ the minimum of $al^d +b/l$ over real numbers $l$ is achieved at $l_0= (b/ad)^{1/(d+1)}$. Taking $l=\lfloor l_0 \rfloor$ when $l_0 \ge 1$ and using that $l_0/2 \le l \le l_0$ yields
\bea \label{eq:min.in.l}
\min_{l \in \mathbb{N}_1} \left(al^d+\frac{b}{l} \right)\le a\left(\frac{b}{ad}\right)^{\frac{d}{d+1}} + 2b\left(\frac{ad}{b}\right)^{\frac{1}{d+1}}
= a^{\frac{1}{d+1}}b^{\frac{d}{d+1}} \left(\frac{1}{d^{\frac{d}{d+1}}}+2 d^{\frac{1}{d+1}}\right).
\ena

\noindent {\textbf{Proof of Theorem \ref{1stappthm}:}}  
By second order stationarity, it suffices to prove the case ${\bf k}={\bf 1}$.
Let $n \ge 2, B_{{\bf 1}}^n$ the block of size $n^d$
as given in \eqref{blockdef}, and $W_{{\bf 1}}^n$ the standardized sum over that block, as in \eqref{wn2}. For any $1 \le l \le n$ write $n=(m-1)l+r, 1 \le r \le l$, and decompose $W_{{\bf 1}}^n$ as the sum of $\xi_{\bf i}/\sqrt{n^dA_n}$ over ${\bf i} \in [m]$, as in \eqref{xidef}.

We apply Theorem \ref{stein:fkg}, dealing with the two terms on the right hand side of \eqref{stein:fkgbound}. For the first term, using $|X_{{\bf j}}| \le K$, the definition \eqref{xidef} of $\xi_{{\bf i}}$, and the fact that the side lengths of all blocks $D_{\bf i}^l$ are at most $l$, we have
\beas 
\left|\frac{\xi_{{\bf i}}}{\sqrt{n^dA_n}} \right| \le B \text{ \ \ with \ } B=\frac{2Kl^d}{\sqrt{n^d A_n}} \text{ \ for all \ } {\bf i} \in [m]^d.
\enas

Applying Lemma \ref{lemma:1stapp} for the last term and invoking Theorem \ref{stein:fkg} now yields
\bea \label{bound1stapp}
d_1\big({\cal L}(W_{{\bf 1}}),{\cal L}(Z)\big) &\le& \frac{10Kl^d}{n^{d/2}A_n^{1/2}}
                     +\frac{5.2\kappa_0 \gamma_{\lambda,d}}{lA_n} .
\ena
Applying the bound \eqref{eq:min.in.l} to the last expression with  $l_0=C_{\lambda,\kappa_0,d}^{-1/(d+1)}n^{d/(2d+2)}$ and $C_{\lambda,\kappa_0,d}$ as in \eqref{cdef} and plugging in $l=\lfloor l_0 \rfloor$ back to the right hand side of \eqref{bound1stapp} yields the result. It is easy to check that $1 \le l_0 \le n$ for $ n \ge \max \left\{ C_{\lambda,\kappa_0,d}^{2/d},C_{\lambda,\kappa_0,d}^{-2/(d+2)} \right\}$.

\bbox

To prove Theorem \ref{1stappthm2}, we apply Theorem \ref{stein:fkg2} and use the same techniques as in Theorem \ref{1stappthm}. We remind the reader that for this result we do not explicitly compute the constants, but index them by the parameters on which they depend.

\noindent {\textbf{Proof of Theorem \ref{1stappthm2}:}}  First we prove the claims that, $\Sigma$ is invertible and $|\Sigma^{-1}|_{\infty}$ is bounded by \eqref{eq:infinity.norm.bound.inverse.lem:Gershgorin} when $n > (p-1)\kappa_0\nu_{\lambda}^de^{-\lambda}/A_n$. Applying Lemma \ref{lemma:1stapp2}, for all $q \in [p]$  we have 
\beas
\Sigma_{q,q} - \sum_{1\le s\le p, s\ne q} \left|\Sigma_{q,s}\right| \ge n^{d-1}(nA_n -(p-1) \kappa_0 \nu_\lambda^d e^{-\lambda}) > 0
\text{ \ if \ } n > \frac{(p-1)\kappa_0\nu_{\lambda}^de^{-\lambda}}{A_n}.
\enas
which implies that $\Sigma$ is a strictly diagonally dominant matrix, and is therefore invertible by the Gershgorin circle theorem, see for instance Theorem 15.10 of \cite{BR14}. Another claim in \eqref{eq:infinity.norm.bound.inverse.lem:Gershgorin} follows from \cite{AN63}, where it is shown that the bound \eqref{eq:infinity.norm.bound.inverse.lem:Gershgorin} holds for the norm $||C||_{\infty} = \max_i \sum_{j =1}^p |c_{ij}|$, which dominates $|C|_{\infty}$.

Next we proceed as in the one dimensional case. For $n \ge 2$, and $1 \le l \le n$ we write $n=(m-1)l+r$ with 
$m \ge 1$ and $1 \le r \le l$, and decompose $S_{{\bf k}_q}^n-\e S_{{\bf k}_q}^n$ for $q\in [p]$ as the sum over ${\bf i} \in [m]^d$ of the variables $\xi_{{\bf i},{\bf k}_q}$ given in \eqref{xidef}.

Applying Theorem \ref{stein:fkg2}, we handle the three terms on the right hand side of (\ref{steinbound:fkg2}). Using the definition \eqref{xidef} of $\xi_{{\bf i},{\bf k}_q}$ and that $|X_{\bf t}|\le K$, we have
\beas 
\left|\xi_{{\bf i},{\bf k}_q}\right| \le  B \text{ \ where \ } B=2Kl^d \text{ \ for all \ } {\bf i} \in [m]^d, q \in [p]
\enas
and thus using $|\Sigma^{-1/2}|_\infty = n^{-d/2}\psi_{n}$ and $\Sigma_{jj}=n^dA_n$, we may bound the first term as
\bea \label{multi1st}
\frac{5}{6}p^3Bn^{-3d/2}\psi_{n}^3\sum_{q=1}^p\Sigma_{q,q} \le \frac{C_{p,K}l^dA_n\psi_{n}^3}{n^{d/2}}.
\ena

For the second term, by Lemma \eqref{lemma:1stapp} we have
\begin{multline} \label{beg}
-\left(\frac{3}{2}p^3Bn^{-3d/2}\psi_{n}^3+p^2n^{-d}\psi_{n}^2\right) \sum_{q=1}^p \sum_{{\bf i},{\bf j} \in [m]^d,{\bf i} \neq {\bf j}} \e \left(\xi_{{\bf i},{\bf k}_q}\xi_{{\bf j},{\bf k}_q}\right) \\
= -\left(\frac{C_{p,K}l^d\psi_{n}^3}{n^{d/2}}+C_p\psi_{n}^2\right) \sum_{q=1}^p \sum_{{\bf i},{\bf j} \in [m]^d,{\bf i} \neq {\bf j}} \e \left(\frac{\xi_{{\bf i},{\bf k}_q}\xi_{{\bf j},{\bf k}_q}}{n^d}\right)
\le \frac{C_{\lambda,\kappa_0,p,K,d}l^{d-1}\psi_{n}^3}{n^{d/2}} + \frac{C_{\lambda,\kappa_0,p,d}\psi_{n}^2}{l}.
\end{multline}

Next, invoking Lemma \ref{lemma:1stapp2} and assumption \eqref{eq:kqks.separated.alphan}
we have
\beas 
-\Sigma_{q,s} = -\Cov \big(S_{{\bf k}_q}^n,S_{{\bf k}_s}^n \big)  \le \kappa_0\nu_\lambda^d e^{-\lambda} n^{d-1} \text{ \ for \ } q \ne s \in [p],
\enas
and hence we may bound the last term as
\bea \label{dest2}
-\left( \frac{2}{3}p^3 B n^{-3d/2}\psi_{n}^3 + p^2 n^{-d}\psi_{n}^2  \right) \sum_{q,s \in [p],q\ne s} \Sigma_{q,s} 
&\le& \left( \frac{C_{p,K}l^d\psi_{n}^3 }{n^{3d/2}} + \frac{C_p\psi_{n}^2}{n^d}\right)\kappa_0\nu_\lambda^d e^{-\lambda} n^{d-1} \nn \\ 
&\le& \frac{C_{\lambda,\kappa_0d,p,K} l^d\psi_{n}^3 }{n^{(d+2)/2}} + \frac{C_{\lambda,\kappa_0,d,p} \psi_{n}^2}{n}.
\ena

By Theorem \ref{stein:fkg2} and \eqref{multi1st}-\eqref{dest2}, we have
\beas 
d_{\mathcal{H}_{3,\infty,p}}\big(\mathcal{L}(\Sigma^{-1/2}\mathbf{W}),\mathcal{L}(\mathbf{Z})\big) &\le& C_{\lambda,\kappa_0,K,p,d}\left( \frac{l^{d-1}\psi_{n}^3}{n^{d/2}} + \frac{ l^d\psi_{n}^3 }{n^{(d+2)/2}} + \frac{ \psi_{n}^2}{n}+ \frac{l^dA_n\psi_{n}^3}{n^{d/2}} + \frac{\psi_{n}^2}{l}\right).
\enas
Since the first three terms do not contribute the rate, applying \eqref{eq:min.in.l} to the last two terms in the parentheses, we obtain 
\beas
l_0 = \left(\frac{1}
{d\psi_n A_n }\right)
^{\frac{1}{d+1}}n^{\frac{d}{2d+2}},
\enas
which satisfies $1 \le l_0 \le n$ for the range of $n$ given in \eqref{eq:thm2rangen}, and pluging in $l=\lfloor l_0 \rfloor$ back to the bound yields the result. 

\bbox

\section*{Acknowledgements} 
The author would like to thank the anonymous referees for extremely helpful comments that lead to an improvement of the bounds appeared in the paper.

\bigskip

\def\cprime{$'$}

\end{document}